\documentclass[12pt]{amsart}

\usepackage{amssymb}
\usepackage{graphicx}
\usepackage{enumerate}
\usepackage{combelow}
\usepackage[letterpaper,margin=1in,ignoreall]{geometry}

\usepackage{hyperref}
\usepackage{url}
\usepackage{breakurl}

\newtheorem*{thm*}{Theorem}
\newtheorem{thm}{Theorem}

\theoremstyle{remark}

\theoremstyle{definition}

\renewcommand{\O}{\mathcal{O}}
\newcommand{\Z}{\mathbb{Z}}

\newcommand{\R}{\mathbb{R}}
\newcommand{\C}{\mathbb{C}}

\newcommand{\J}{\mathcal{J}}

\newcommand{\One}{\mathbf{1}}

\DeclareMathOperator{\lct}{lct}
\DeclareMathOperator{\Newt}{Newt}
\DeclareMathOperator{\Int}{Int}
\DeclareMathOperator{\monom}{monom}

\newcommand{\sw}[1]{\textsc{#1}} 

\title{Software for multiplier ideals}

\author{Zach Teitler}
\address{Department of Mathematics \\
1910 University Drive \\
Boise State University \\
Boise, ID 83725-1555 \\
USA}
\email{zteitler@boisestate.edu}

\date{\today}



\begin{document}

\bibliographystyle{amsplain}       

\begin{abstract}
We describe a new software package for computing multiplier ideals
in certain cases, including monomial ideals, monomial curves, generic determinantal ideals,
and hyperplane arrangements.
In these cases we take advantage of combinatorial formulas for multiplier ideals
given by results of Howald, Thompson, and Johnson.
The package uses \sw{Normaliz}~\cite{MR2659215,MR2881130}.
It is available as a library for \sw{Macaulay2}~\cite{M2}.
\end{abstract}

\maketitle


\section{Introduction}

Multiplier ideals have been applied to a number of problems in algebraic geometry in recent years,
most spectacularly in recent major advances in the minimal model program \cite{MR2359343,MR2601039}
that built on earlier work showing the deformation invariance of plurigenera \cite{siu}.
Other applications include several results on singularities and linear series
(see \cite{pag2} and \cite{MR2275611}),
a bound for symbolic powers (see \cite{els:symbolic-powers}),
and applications to algebraic statistics
(see for example \cite{MR2554932}, \cite{MR2877601}, \cite[Chapter 5]{MR2723140}).
New applications of multiplier ideals continue to emerge in topics such as Chow stability \cite{MR2389549}
and singularities in generic liaison \cite{Niu:2012fk}.
With broad and growing interest in multiplier ideals, it is increasingly valuable to compute examples.

Here is a definition of multiplier ideals in terms of resolution of singularities.
Suppose $X$ is a smooth complex variety
(which we may assume is affine, or even just $\C^n$, since we are primarily interested in local issues),
$I \subset \O_X$ is a nonzero ideal sheaf,
and $\mu : Y \to X$ is a log resolution of $I$,
so that the total transform $I \O_Y$ defines a divisor $F$ with simple normal crossings support,
$F = \sum a_i E_i$, where the $E_i$ are distinct reduced components of $F$.
Then for each real number $c \geq 0$, the $c$'th multiplier ideal is defined by
$  \J(I^c) = \mu_* \O_Y ( K_{Y/X} - \lfloor c \cdot F \rfloor )$
where $K_{Y/X}$ is the relative canonical divisor of $Y$ over $X$, defined locally by the vanishing of the
determinant of the Jacobian $d\mu$, and $\lfloor c \cdot F \rfloor$ denotes the component-wise round-down
of the $\R$-divisor $c \cdot F$, given by $\lfloor c \cdot F \rfloor = \sum \lfloor c a_i \rfloor E_i$.
There are a number of equivalent characterizations of multiplier ideals, in terms of
jet spaces~\cite{MR2081163},
$D$-modules~\cite{MR2123230},
test ideals for tight closure~\cite{MR1808611},
and local integrability~\cite{MR0466635,MR512213,MR1078269}.
For a complete account of multiplier ideals emphasizing the approach via resolution of singularities, see \cite{pag2}.

In theory it is algorithmic to compute multiplier ideals by computing a resolution of singularities of $I$ followed by a sheaf pushforward.
In practice it is more difficult, see \cite{Fruhbis-Kruger:2013fk}.

Shibuta's algorithm for computing Bernstein-Sato polynomials and multiplier ideals
via Gr\"obner basis methods in Weyl algebras \cite{Shibuta20112829}
(implemented by Shibuta in \sw{Risa/Asir})
was refined and implemented in the \sw{Dmodules} library for \sw{Macaulay2}
by Berkesch and Leykin \cite{DBLP:conf/issac/BerkeschL10}.
The \sw{Dmodules} library can compute multiplier ideals and jumping numbers of arbitrary ideals,
but due to the difficulty of the computations, can only handle modestly sized examples.

We describe a new software package named \sw{MultiplierIdeals}
that computes multiplier ideals of special ideals
including monomial ideals, ideals of monomial curves, generic determinantal ideals,
and hyperplane arrangements
via combinatorial methods,
using the \sw{Normaliz} software and interface to \sw{Macaulay2} by Bruns, et al
\cite{MR2659215,MR2881130}.
The combinatorial methods allow computations of somewhat larger examples than can be handled
by general methods.

Our package also computes certain quantities associated to multiplier ideals,
the log canonical thresholds and jumping numbers.
Because of the round-down operation, $\J(I^{c+\epsilon}) = \J(I^c)$ for sufficiently small $\epsilon > 0$.
A real number $c \geq 0$ is a jumping number of $I$ if $\J(I^c) \neq \J(I^{c-\epsilon})$ for all $\epsilon > 0$.
Every jumping number is in fact rational.
The smallest strictly positive jumping number is called the log canonical threshold of $I$, denoted $\lct(I)$.
It turns out that $\J(I^0) = (1)$ is the trivial ideal, so $\lct(I)$ is the supremum of $c$
such that $\J(I^c) = (1)$; equivalently, $\lct(I)$ is the first value of $c$ such that $\J(I^c) \neq (1)$.

The portion dealing with monomial ideals was written first
and distributed as a package named \sw{MonomialMultiplierIdeals}.
The portion of the package dealing with monomial curves was written
by Claudiu Raicu, Bart Snapp, and the author
at the 2011 IMA Special Workshop on Macaulay2,
and distributed as a package named \sw{SpaceCurvesMultiplierIdeals}.
The portion dealing with hyperplane arrangements is based on
code written by Denham and Smith
for the \sw{HyperplaneArrangements} package \cite{M2-HyperplaneArrangements}.
These portions were all integrated into the present package,
and computations with generic determinantal ideals added,
at the 2012 Macaulay2 Workshop at Wake Forest.

The \sw{MultiplierIdeals} package is available from the author's web site at
\url{http://math.boisestate.edu/~zteitler/math/MultiplierIdealsSoftware.php}.


\section{Monomial ideals}

For a monomial ideal $I \subset \C[x_1,\dotsc,x_n]$, let $\monom(I) \subset \Z_{\geq 0}^n$ be the set of exponent vectors of monomials in $I$.
The Newton polyhedron $\Newt(I)$ is the convex hull of $\monom(I)$.
Let $\One=(1,\dotsc,1) \in \R^n$.

Multiplier ideals of monomial ideals are described by the following theorem of Howald.
\begin{thm}[\cite{howald:monomial}]
The multiplier ideal $\J(I^c)$ is the monomial ideal containing $x^v$ if and only if
$v + \One \in \Int(c \cdot \Newt(I))$.
Here $\Int$ denotes the topological interior of $c \cdot \Newt(I)$ relative to the nonnegative orthant,
that is, as a subset of $(\R_{\geq0})^n$.
\end{thm}

In other words, the multiplier ideal is the quotient ideal
\[
  \J(I^c) = ( x^v : v \in \Int(c \cdot \Newt(I)) ) : x^\One .
\]
The Newton polyhedron $\Newt(I)$ is defined by
a system of inequalities $Av \geq b$ where $A$ is an $r \times n$ matrix,
$b$ is a vector,
and $\geq$ is the partial order of entrywise comparison,
where $a \geq b$ if and only if $a_i \geq b_i$, $1 \leq i \leq r$.
Then $c \cdot \Newt(I)$ is defined by $Av \geq cb$.
The interior $\Int(c \cdot \Newt(I))$ is the solution of the system of inequalities given by
\[
  A_i v > c b_i \quad \text{if $b_i \neq 0$},
  \qquad\qquad
  A_i v \geq c b_i = 0 \quad \text{if $b_i = 0$}.
\]

Since $\Newt(I)$ is a rational polyhedron, we can (and do) take the $A$ and $b$ to have integer entries.
Furthermore since $I$ is an ideal the entries of $A$ and $b$ are nonnegative.
In practice, it is sufficient to compute $\J(I^c)$ for rational $c = p/q$, and this can be done as follows.
To find the integer vectors $v$ lying in the topological interior of the solution region to $Av \geq cb$, equivalently $qAv \geq pb$,
we add $1$ to the nonzero entries of $p b$, yielding a vector $b'$ with entries $b'_i = p b_i + 1$ if $b_i \neq 0$, otherwise $b'_i = p b_i = b_i = 0$.
Then the multiplier ideal $\J(I^c)$ is the quotient $( x^v : qAv \geq b' ) : x^\One$.

The software \sw{Normaliz} can compute the defining inequalities $Av \geq b$ of $\Newt(I)$
and the solutions to the modified system $qAv \geq b'$;
\sw{Macaulay2} can compute the ideal quotient by $x^\One$,
giving the multiplier ideal.

\begin{verbatim}
$ M2
Macaulay2, version 1.6
i1 : loadPackage "MultiplierIdeals";
i2 : R = QQ[x,y,z,w];
i3 : I = monomialIdeal(x*y, x*z, y*z, y*w, z*w^2);
o3 : MonomialIdeal of R
i4 : logCanonicalThreshold(I)
o4 = 2
i5 : multiplierIdeal(I,7/3)
                     2
o5 = ideal (y, z*w, z , x*z)
i6 : toString jumpingNumbers(I)
o6 = {{2, 7/3, 5/2, 8/3, 3, 10/3, 7/2, 11/3, 4}, {ideal(z,y),
     ideal(y,z*w,z^2,x*z), ideal(z*w,y*w,y*z,x*z,y^2,x*y),
     ideal(y*w,y*z,x*z,y^2,x*y,z*w^2,z^2*w),
     ideal(y*z*w,y^2*w,y*z^2,x*z^2,y^2*z,x*y*z,x*y^2,z^2*w^2), ...
\end{verbatim}
The \verb+jumpingNumbers+ command produces a list with two elements:
\begin{enumerate}
\item A list of the jumping numbers of $I$ in the interval $(0,k(I)]$,
where $k(I)$ is the analytic spread of $I$.
A different interval may be specified as an optional argument.
\item A list of the multiplier ideals at the jumping numbers.
(The list is truncated in the above example.)
\end{enumerate}
Thus the output of the last command says that this ideal $I$ has jumping numbers $2, 7/3,\dotsc$,
and gives the corresponding multiplier ideals:
$\J(I^2) = (z,y)$, $\J(I^{7/3}) = (y,zw,z^2,xz)$, and so on.

Multiplier ideals and jumping numbers for $c > k(I)$ are given by Skoda's theorem
\cite[Thm.~9.6.21]{pag2}.

In the above example, the log canonical threshold, single multiplier ideal $\J(I^{7/3})$,
and list of $9$ jumping numbers and multiplier ideals were each computed
in a fraction of a second
on a 2012 MacBook with dual-core 64-bit 2.9 GHz CPU and 8GB RAM.
By way of comparison, the \sw{Dmodules} package takes about 42 seconds to compute the log canonical
threshold on the same machine, and about 84 seconds to compute $\J(I^{7/3})$.
Of course, it is in no way a criticism of the \sw{Dmodules} package that it takes longer;
it is using a general method that works for any ideal,
while our package uses a very special reduction to a combinatorial algorithm for monomial ideals.
The comparison is intended only to illustrate the advantages
of using special algorithms where available.

For monomial ideals, extra information is available:
for any monomial $x^v$, the package computes the threshold value
$\min \{ c : x^v \notin \J(I^c) \}$, and the list of facets of the Newton polyhedron
that impose the nonmembership.
\begin{verbatim}
i7 : toString logCanonicalThreshold(I,z^2*w)
o7 = (3,matrix {{2, 2, 1, 1, -3}, {2, 2, 0, 1, -2}})
\end{verbatim}
This output means that $z^2w \notin \J(I^3)$ but $z^2w \in \J(I^c)$ for $c < 3$.
That is, for the exponent vector $v = (0,0,2,1)$, $v + \One$ lies on the boundary of $3 \cdot \Newt(I)$;
and furthermore it lies on the intersection of two facets,
the ones scaled up from the facets of $\Newt(I)$ defined by
$2x+2y+z+w=3$ and $2x+2y+w=2$.

The log canonical threshold of the ideal $I$ itself is the threshold
value for $1 = x^0$.


\section{Monomial curves}

An affine monomial curve is one parametrized by $t \mapsto (t^{a_1},\dotsc,t^{a_n})$.
We may and do assume $1 \leq a_1 \leq \dotsc \leq a_n$ and $\gcd(a_1,\dotsc,a_n) = 1$.
For convenience we denote this curve $C(a_1,\dotsc,a_n)$.
It has a singularity at the origin when $a_1 \geq 2$.
The defining ideal is the kernel of the map $\C[x_1,\dotsc,x_n] \to \C[t]$
given by $x_i \mapsto t^{a_i}$, a binomial ideal.

The multiplier ideals of affine monomial curves in dimension $n=3$ have been found
by Howard Thompson \cite{Thompson:2010uq},
using the combinatorial description of the resolution of singularities of a binomial ideal
given in \cite{MR1892938}.
This yields a combinatorial formula in terms of the vector $(a,b,c)$ of exponents
appearing in the parametrization $t \mapsto (t^a,t^b,t^c)$.
Our software package implements Thompson's result, again calling on \sw{Normaliz}
to find generators for the semigroup of integer solutions to certain linear inequalities.

\begin{verbatim}
$ M2
i1 : loadPackage "MultiplierIdeals";
i2 : R = QQ[x,y,z]; S = QQ[t];
i3 : I = kernel map(S,R,{t^3,t^4,t^5})
             2         2     2   3
o3 = ideal (y  - x*z, x y - z , x  - y*z)
\end{verbatim}
To compute the multiplier ideals and log canonical threshold of $I$,
we input the list of exponents in the parameterization:
\begin{verbatim}
i4 : toString logCanonicalThreshold(R,{3,4,5})
o4 = 13/9
\end{verbatim}
We compare with the log canonical threshold of the generic initial ideal (gin) of $I$, a monomial ideal:
\begin{verbatim}
i5 : loadPackage "GenericInitialIdeal";
i6 : gin I
             2     2   3
o6 = ideal (x , x*y , y )
i7 : toString logCanonicalThreshold monomialIdeal gin I
o7 = 5/6
\end{verbatim}
We can compute the generic initial ideal of the same ideal with respect to the lex order,
obtaining a different initial ideal with a different log canonical threshold:
\begin{verbatim}
i8 : lexgin I
             2        5     4
o8 = ideal (x , x*y, y , x*z )
i9 : toString logCanonicalThreshold monomialIdeal lexgin I
o9 = 21/20
\end{verbatim}

It would be interesting to investigate the relationships between
properties of the numerical semigroup generated by $(a,b,c)$,
the multiplier ideals of the monomial curve $C(a,b,c)$,
and the multiplier ideals of the generic initial ideal of $C(a,b,c)$
(and in addition to multiplier ideals, log canonical threshold and jumping numbers).

To indicate how such an investigation might begin,
the following table shows the log canonical thresholds of some monomial curves $C(a,b,c)$
and the log canonical thresholds of the generic initial ideals of the same curves.
It also shows the Frobenius number of the semigroup generated by $(a,b,c)$,
the genus of the semigroup,
and and whether the semigroup is symmetric, as computed by the \sw{Numericalsgps} package
for \sw{GAP} \cite{GAP-NumericalSemigps,GAP4}.
\[
\begin{array}{l | lllll}
(a,b,c) & \lct & \lct \operatorname{gin} & \text{Frobenius} & \text{genus} & \text{Symmetric} \\ \hline
(3, 4, 5)   & 13/9  & 5/6  & 2 & 2 & \\
(3, 5, 7)   & 17/12 & 3/4  & 4 & 3 & \\
(3, 7, 8)   & 5/4   & 2/3  & 5 & 4 & \\
(3, 7, 11)  & 25/18 & 2/3  & 8 & 5 & \\
(3, 8, 10)  & 5/4   & 7/12 & 7 & 5 & \\
(3, 10, 11) & 7/6   & 7/13 & 8 & 6 & \\
(4, 5, 6)   & 17/12 & 3/4  & 7 & 4 & \text{yes} \\
(4, 6, 7)   & 4/3   & 2/3  & 9 & 5 & \text{yes} \\
\end{array}
\]


\section{Generic determinantal ideals}

Let $X = (x_{i,j})_{\substack{1 \leq i \leq m \\ 1 \leq j \leq n}}$ be an $m \times n$ generic matrix,
meaning one whose entries are independent variables.
Let $I_r(X)$ be the ideal generated by $r \times r$ minors of $X$.
The multiplier ideals of $I_r(X)$ have been found by Amanda Johnson in her
2003 dissertation \cite{amanda}.
\begin{thm}[\cite{amanda}]
With $X$, $m$, $n$, $r$ as above, over an algebraically closed field,
the multiplier ideals are given by the following intersection of symbolic powers
of determinantal ideals:
\[
  \J(I_r(X)^c) = \bigcap_{i=1}^{r} I_i(X)^{(\lfloor c(r+1-i) \rfloor + 1 - (n-i+1)(m-i+1) )} .
\]
\end{thm}
Recall that symbolic powers of generic determinantal ideals may be expressed as
\[
  I_r(X)^{(a)} = \sum_{\kappa_1+\dotsb+\kappa_s=a} \prod_{i=1}^s I_{r-1+\kappa_i}(X) ,
\]
the sum taken over partitions of $a$.
See \cite[Theorem 10.4]{MR953963}.

We may compute multiplier ideals of determinantal ideals in our software
by giving the matrix $X$ and the size of minors.
Here we examine multiplier ideals of the size $2$ minors and size $3$ minors
of a $4 \times 5$ generic matrix.
\begin{verbatim}
i1 : loadPackage "MultiplierIdeals";
i2 : R = QQ[x_1..x_20];
i3 : X = genericMatrix(R,4,5); -- a 4x5 generic matrix
i4 : logCanonicalThreshold(X,2) -- lct of the ideal of 2x2 minors
o4 = 10
i5 : multiplierIdeal(X,2,10) == minors(1,X) -- J(I^10) where I = 2x2 minors
o5 = true
i6 : multiplierIdeal(X,2,11) == (minors(1,X))^3 -- J(I^11)
o6 = true
i7 : logCanonicalThreshold(X,3) -- lct of the ideal of 3x3 minors
o7 = 6
i8 : multiplierIdeal(X,3,6) == minors(3,X) -- J(I^6) where I = 3x3 minors
o8 = true
\end{verbatim}


\section{Hyperplane arrangements}

A formula for multiplier ideals of hyperplane arrangements was found by
Musta\cb{t}\u{a} \cite{mustata:hyperplane-arrangements}
and simplified in \cite{teitler:hyperplane-arrangements}.
The \sw{HyperplaneArrangements} package \cite{M2-HyperplaneArrangements}
uses these results to compute multiplier ideals and log canonical thresholds of
hyperplane arrangements.
To this we add the ability to compute jumping numbers
and other minor modifications.
I thank the authors of \sw{HyperplaneArrangements},
Graham Denham and Gregory G.\ Smith,
for their permission to copy and modify their package's source code.

The following is Example~6.3 of \cite{DBLP:conf/issac/BerkeschL10}.
\begin{verbatim}
i1 : loadPackage "MultiplierIdeals";
i2 : R = QQ[x,y,z];
i3 : ff = toList factor ( (x^2-y^2)*(x^2-z^2)*(y^2-z^2)*z ) / first;
i4 : A = arrangement ff;
i5 : toString jumpingNumbers(A,IntervalType=>"ClosedOpen")
o5 = {{3/7, 4/7, 2/3, 6/7}, {ideal(z,y,x), ideal(z^2,y*z,x*z,y^2,x*y,x^2),
  ideal(y^2*z-z^3,x^2*z-z^3,x*y^2-x*z^2,x^2*y-y*z^2), ...
\end{verbatim}


\section*{Acknowledgements}

I am very grateful to Claudiu Raicu and Bart Snapp for their critical contributions
to the package and for a number of very helpful comments
about this paper and the software package itself.
I would also like to thank Howard Thompson for sharing his work in progress
and for numerous helpful conversations,
Graham Denham and Gregory G.\ Smith,
the organizers of the 2011 IMA Special Workshop on Macaulay2,
and the organizers of the 2012 Macaulay2 Workshop at Wake Forest.


\bigskip

\renewcommand{\MR}[1]{}
\bibliography{/Users/zteitler/Math/research/biblio}

\providecommand{\bysame}{\leavevmode\hbox to3em{\hrulefill}\thinspace}
\providecommand{\MR}{\relax\ifhmode\unskip\space\fi MR }
\providecommand{\MRhref}[2]{%
  \href{http://www.ams.org/mathscinet-getitem?mr=#1}{#2}
}
\providecommand{\href}[2]{#2}
\begin{thebibliography}{10}

\bibitem{DBLP:conf/issac/BerkeschL10}
Christine Berkesch and Anton Leykin, \emph{Algorithms for {B}ernstein--{S}ato
  polynomials and multiplier ideals}, ISSAC (Wolfram Koepf, ed.), ACM, 2010,
  pp.~99--106.

\bibitem{MR2601039}
Caucher Birkar, Paolo Cascini, Christopher~D. Hacon, and James
  M{\textsuperscript{c}}Kernan, \emph{Existence of minimal models for varieties
  of log general type}, J. Amer. Math. Soc. \textbf{23} (2010), no.~2,
  405--468. \MR{2601039}

\bibitem{MR2659215}
Winfried Bruns and Bogdan Ichim, \emph{Normaliz: algorithms for affine monoids
  and rational cones}, J. Algebra \textbf{324} (2010), no.~5, 1098--1113.
  \MR{2659215}

\bibitem{MR2881130}
Winfried Bruns and Gesa K{\"a}mpf, \emph{A {M}acaulay2 interface for
  {N}ormaliz}, J. Softw. Algebra Geom. \textbf{2} (2010), 15--19. \MR{2881130}

\bibitem{MR953963}
Winfried Bruns and Udo Vetter, \emph{Determinantal rings}, Lecture Notes in
  Mathematics, vol. 1327, Springer-Verlag, Berlin, 1988. \MR{953963
  (89i:13001)}

\bibitem{MR2123230}
Nero Budur and Morihiko Saito, \emph{Multiplier ideals, {$V$}-filtration, and
  spectrum}, J. Algebraic Geom. \textbf{14} (2005), no.~2, 269--282.
  \MR{2123230 (2006g:14012)}

\bibitem{GAP-NumericalSemigps}
Manuel Delgado and Pedro~A. Garc{\'\i}a-S{\'a}nchez, \emph{{Numericalsgps -- a
  GAP package, Version 0.971}},
  \url{http://www.gap-system.org/Packages/numericalsgps.html}, 2011.

\bibitem{M2-HyperplaneArrangements}
Graham Denham and Gregory~G. Smith, \emph{{HyperplaneArrangements -- a package
  for Macaulay2, version 0.9}},
  \url{http://www.math.uiuc.edu/Macaulay2/doc/Macaulay2-1.6/share/doc/Macaulay2/HyperplaneArrangements/html/},
  2011.

\bibitem{MR2723140}
Mathias Drton, Bernd Sturmfels, and Seth Sullivant, \emph{Lectures on algebraic
  statistics}, Oberwolfach Seminars, vol.~39, Birkh\"auser Verlag, Basel, 2009.
  \MR{2723140 (2012d:62004)}

\bibitem{MR2081163}
Lawrence Ein, Robert Lazarsfeld, and Mircea Musta{\c{t}}{\v{a}}, \emph{Contact
  loci in arc spaces}, Compos. Math. \textbf{140} (2004), no.~5, 1229--1244.
  \MR{2081163 (2005f:14006)}

\bibitem{els:symbolic-powers}
Lawrence Ein, Robert Lazarsfeld, and Karen~E. Smith, \emph{Uniform bounds and
  symbolic powers on smooth varieties}, Invent. Math. \textbf{144} (2001),
  no.~2, 241--252. \MR{MR1826369 (2002b:13001)}

\bibitem{MR2275611}
Lawrence Ein and Mircea Musta{\c{t}}{\u{a}}, \emph{Invariants of singularities
  of pairs}, International {C}ongress of {M}athematicians. {V}ol. {II}, Eur.
  Math. Soc., Z{\"u}rich, 2006, pp.~583--602. \MR{MR2275611 (2007m:14050)}

\bibitem{Fruhbis-Kruger:2013fk}
Anne Fr{{\"u}}hbis-Kr{{\"u}}ger, \emph{Desingularization in computational
  applications and experiments},
  \href{http://arxiv.org/abs/1301.3709}{\nolinkurl{arXiv:1301.3709}} [math.AG],
  Jan 2013.

\bibitem{GAP4}
The GAP~Group, \emph{{GAP -- Groups, Algorithms, and Programming, Version
  4.6.4}}, 2013.

\bibitem{MR1892938}
Pedro~Daniel Gonz{\'a}lez~P{\'e}rez and Bernard Teissier, \emph{Embedded
  resolutions of non necessarily normal affine toric varieties}, C. R. Math.
  Acad. Sci. Paris \textbf{334} (2002), no.~5, 379--382. \MR{1892938
  (2003b:14019)}

\bibitem{M2}
Daniel~R. Grayson and Michael~E. Stillman, \emph{Macaulay2, a software system
  for research in algebraic geometry}, Available at
  \href{http://www.math.uiuc.edu/Macaulay2/}%
  {http://www.math.uiuc.edu/Macaulay2/}.

\bibitem{MR2359343}
Christopher~D. Hacon and James M{\textsuperscript{c}}Kernan, \emph{Extension
  theorems and the existence of flips}, Flips for 3-folds and 4-folds, Oxford
  Lecture Ser. Math. Appl., vol.~35, Oxford Univ. Press, Oxford, 2007,
  pp.~76--110. \MR{MR2359343}

\bibitem{howald:monomial}
J.~A. Howald, \emph{Multiplier ideals of monomial ideals}, Trans. Amer. Math.
  Soc. \textbf{353} (2001), no.~7, 2665--2671 (electronic). \MR{MR1828466
  (2002b:14061)}

\bibitem{amanda}
Amanda~A. Johnson, \emph{Multiplier ideals of determinantal ideals}, Ph.D.
  thesis, U. Michigan, 2003.

\bibitem{MR0466635}
J.~J. Kohn, \emph{Sufficient conditions for subellipticity on weakly
  pseudo-convex domains}, Proc. Nat. Acad. Sci. U.S.A. \textbf{74} (1977),
  no.~6, 2214--2216. \MR{0466635 (57 \#6512)}

\bibitem{MR512213}
\bysame, \emph{Subellipticity of the {$\bar \partial $}-{N}eumann problem on
  pseudo-convex domains: sufficient conditions}, Acta Math. \textbf{142}
  (1979), no.~1-2, 79--122. \MR{512213 (80d:32020)}

\bibitem{pag2}
Robert Lazarsfeld, \emph{Positivity in algebraic geometry. {II}}, Ergebnisse
  der Mathematik., vol.~49, Springer-Verlag, Berlin, 2004, Positivity for
  vector bundles, and multiplier ideals. \MR{MR2095472}

\bibitem{MR2389549}
Yongnam Lee, \emph{Chow stability criterion in terms of log canonical
  threshold}, J. Korean Math. Soc. \textbf{45} (2008), no.~2, 467--477.
  \MR{2389549 (2009b:14091)}

\bibitem{mustata:hyperplane-arrangements}
Mircea Musta{\c{t}}{\v{a}}, \emph{Multiplier ideals of hyperplane
  arrangements}, Trans. Amer. Math. Soc. \textbf{358} (2006), 5015--5023.

\bibitem{MR1078269}
Alan~Michael Nadel, \emph{Multiplier ideal sheaves and {K}{\"a}hler-{E}instein
  metrics of positive scalar curvature}, Ann. of Math. (2) \textbf{132} (1990),
  no.~3, 549--596. \MR{1078269 (92d:32038)}

\bibitem{Niu:2012fk}
Wenbo Niu, \emph{Singularities of generic linkage of algebraic varieties},
  \href{http://arxiv.org/abs/1207.1082}{\nolinkurl{arXiv:1207.1082}} [math.AG],
  Jul 2012.

\bibitem{Shibuta20112829}
Takafumi Shibuta, \emph{Algorithms for computing multiplier ideals}, Journal of
  Pure and Applied Algebra \textbf{215} (2011), no.~12, 2829--2842.

\bibitem{siu}
Yum-Tong Siu, \emph{Invariance of plurigenera}, Invent. Math. \textbf{134}
  (1998), no.~3, 661--673. \MR{MR1660941 (99i:32035)}

\bibitem{MR1808611}
Karen~E. Smith, \emph{The multiplier ideal is a universal test ideal}, Comm.
  Algebra \textbf{28} (2000), no.~12, 5915--5929, Special issue in honor of
  Robin Hartshorne. \MR{1808611 (2002d:13008)}

\bibitem{teitler:hyperplane-arrangements}
Zach Teitler, \emph{A note on {M}usta{\c t}{\u a}'s computation of multiplier
  ideals of hyperplane arrangements}, Proc. Amer. Math. Soc. \textbf{136}
  (2008), no.~5, 1575--1579.

\bibitem{Thompson:2010uq}
Howard~M Thompson, \emph{Multiplier ideals of monomial space curves},
  \href{http://arxiv.org/abs/1006.1915}{\nolinkurl{arXiv:1006.1915v4}}
  [math.AG], Jun 2010.

\bibitem{MR2554932}
Sumio Watanabe, \emph{Algebraic geometry and statistical learning theory},
  Cambridge Monographs on Applied and Computational Mathematics, vol.~25,
  Cambridge University Press, Cambridge, 2009. \MR{2554932 (2011g:62185)}

\bibitem{MR2877601}
Piotr Zwiernik, \emph{An asymptotic behaviour of the marginal likelihood for
  general {M}arkov models}, J. Mach. Learn. Res. \textbf{12} (2011),
  3283--3310. \MR{2877601 (2012m:62248)}

\end{thebibliography}

\end{document}